\documentclass{article}

\usepackage{amssymb} 
\usepackage{amsmath} 
\usepackage{latexsym} 

\usepackage{epsfig}

\def\e{\varepsilon}
\def\p{\prime}
\def\w{\wedge}
\def\vect#1{\mbox{\boldmath $#1$}} 
\def\omegasub#1{\mbox{\large $\omega$}\mbox{\small${\mbox{\large ${}$}}_{#1}$}}
\def\nekofrac#1#2{\mbox{\footnotesize$\displaystyle\frac{#1}{#2}$}}
\def\Si{\varSigma}

\newtheorem{theorem}{Theorem}[section]
\newtheorem{lemma}[theorem]{Lemma}
\newtheorem{corollary}[theorem]{Corollary}
\newtheorem{proposition}[theorem]{Proposition}

\newtheorem{definition}[theorem]{Definition}

\newenvironment{remark}{{\noindent \bf Remark:}}{}

\newenvironment{proof}{{\noindent \bfseries\itshape Proof: }}{\hfill$\Box$} 

\begin{document}

\setlength\arraycolsep{1pt} 

\title{A Note on Y-energies of Knots}

\author{Jun O'Hara}

\maketitle

\begin{abstract} {\footnotesize
We study a $1$-form which can be given by a vector in a conformally invariant way. 
We then study conformally invariant functionals associated to a ``Y-diagram'' on the space of knots which are made from the $1$-form. }
\end{abstract}





\section{Introduction} 
The study was motivated to merge the following two: 

\smallskip
In \cite{OH1} we defined the {\em energy} of a knot. 
It can be considered as the normalization of modified electrostatic energy of a charged knot. 
It measures geometric complexity of knots. 
Later on, it was proved to be conformally invariant (\cite{Fr-He-Wa}). 

\smallskip
On the other hand, Lin and Wang gave functionals associated to choord diagrams on the space of knots using a $1$-form which comes from the Gauss integral formula for the linking number (\cite{Li-Wa}). 

\smallskip
The energy of a knot can be considered as the integration of the interaction between a pair of points on a knot. 
Let us proceed to study functionals so that more than two points on a knot are involved in the integrands. 
For this purpose, we study functionals associated to a chord diagram of the shape of ``Y''. 

We first show that a vector in $\mathbb R^3$ can give, in a conformally invariant way, a vector field, or equivalently, a $1$-form on the complement of the point from which the vector starts. 
This can be done by a {\em conformal transportation} of the vector, which is generalization of the notion introduced in \cite{La-OH}. 
We can play a similar game to that in \cite{Li-Wa} by substituting this conformally invariant $1$-form for the $1$-form coming from the Gauss formula for the linking number. 
But unfortunately, it does not work well: 
It turns out that our functionals are either identically $0$ or identicaly $+\infty$ for any knot. 

\medskip
This article serves as an errata to the author's talk at the AMS Sprin Western Sectional Meeting at San Francisco in May 2003. 
He presented four kinds of $Y$-energies; two are confomally invariant and 
the other two not. 
But the first one turns out to be $0$, the second one diverges, and the third one was already given in Lin and Wang's paper, and the forth one can be obtained from the third one by replacing the integrand by its absolute value.

\section{The results by Lin and Wang}\label{Lin-Wang}
Let us start with introducing the study by Lin and Wang \cite{Li-Wa}. 

\subsection{The linking number} 
Let $L=K_1\cup K_2$ be a $2$-component link in $\mathbb{R}^3$ 
with $K_1=f(S^1)$ and $K_2=g(S^1)$. 
Put 
\begin{equation}\label{def_Gauss_map}
\varphi:\mathbb R^3\times\mathbb R^3\setminus \Delta\ni (x,y)
\mapsto \frac{x-y}{|x-y|}\in S^2. 
\end{equation}
Define a map $\varphi_L$ from a torus $K_1\times K_2$ to $S^2$ by $\varphi_L=\varphi|_{K_1\times K_2}$. 
The linking number $\mbox{{\sl Lk}}(K_1, K_2)$ of $K_1$ and $K_2$ is equal to the degree of $\varphi_{{}_L}$. 

Let $\omegasub{S^2}$ be the unit volume form of $S^2$: 
\begin{equation}\label{area_form_s^2}
\omegasub{S^2}
=\frac{\,1\,}{4\pi}\,\frac{x_1dx_2\w dx_3+x_2dx_3\w dx_1+x_3dx_1\w dx_2}
{|x|^3}. 
\end{equation}
Define a $2$-form $\omega$ on $\mathbb R^3\times\mathbb R^3\setminus \Delta$ by $\omega=\varphi^{\ast}\omegasub{S^2}$. 
Then the linking number can be expressed by the Gauss integral: 
\begin{eqnarray}
GI(f,g)&=&\displaystyle \int _{K_1\times K_2}\varphi_{{}_L}{}^{\ast}\omegasub{S^2}=\int _{K_1\times K_2}\omega(x,y) \nonumber\\[2mm]
&=&\displaystyle \frac{\,1\,}{4\pi}\iint_{S^1\times S^1}\frac
{\det(f^{\p}(s),g^{\p}(t),f(s)-g(t))}{|f(s)-g(t)|^3}\,dsdt\,, 
\label{formula_Gauss_integral}
\end{eqnarray}
where $f(s), g(t)$ are considered as column vectors 
and $\times$ denotes the vector product in $\mathbb{R}^3$. 

\subsection{Integration associated with chord diagrams}

The formulae (\ref{def_energy}) and (\ref{formula_Gauss_integral}) can be considered as the integrals of the interactions between a pair of points on a knot or a link. 
We can generalize them by taking into account more complex combinations of points. 

\begin{figure}[htbp]
\begin{center}
\includegraphics[width=.25\linewidth]{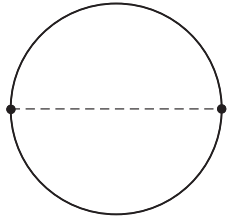}
\caption{The ``$\theta$-graph"}
\label{fig_theta-graph}
\end{center}
\end{figure}

For this purpose, it is natural to use so-called chord diagrams which are used in the study of the Kontsevich integral of the Vassiliev invariant. 
Then (\ref{def_energy}) can be considered to be associated with the ``$\theta$-graph" as illustrated in Figure \ref{fig_theta-graph}, and (\ref{formula_Gauss_integral}) with the ``handcuffs graph". 

\subsection{Gauss integral asociated with chord diagrams}

It is convenient to give the domain of the integral associated with a chord diagram as a subset of a configuration space. 
Let $p_1, \cdots, p_n$ be points on an oriented knot. 
We write $p_1\prec \cdots \prec p_n$\index{$\prec$} if the cyclic order of 
$p_1, \cdots, p_n$ coincides with the orientation of the knot. 
Put 
\begin{equation}\label{def_U_Y}
\begin{array}{rcl}
{\rm Conf}_{3,1}(K,\mathbb{R}^3)&=&\displaystyle \left.\left\{(y_1, y_2, y_3, x)\,\right|
y_1,y_2,y_3\in \!K,\, x\!\in\mathbb{R}^3\setminus K,\, 
y_j\ne y_k \, (j\ne k)\right\},\\[2mm]
U_Y&=&\displaystyle \left.\left\{(y_1, y_2, y_3, x)\in{\rm Conf}_{3,1}(K,\mathbb{R}^3)\,\right|\,
y_1\prec y_2 \prec y_3\right\}\,.
\end{array}\end{equation}

\begin{definition} \rm (\cite{Li-Wa}) 
Let $\omega=\varphi^{\ast}\omegasub{S^2}$ as before. 
Define the $X$-Gauss integral and $Y$-Gauss integral of a knot $K$ by 
\begin{eqnarray}\label{xy}
GI_X(K)&=&\displaystyle \int_{y_1\prec y_2 \prec y_3 \prec y_4}
\omega(y_1, y_3)\wedge\omega(y_2, y_4), \\[2mm]
GI_Y(K)&=&\displaystyle \int_{U_Y}\omega(x, y_1)\wedge\omega(x, y_2)\wedge\omega(x, y_3).
\end{eqnarray}
\end{definition}

\begin{figure}[htbp]
\begin{center}
\begin{minipage}{.45\linewidth}
\begin{center}
\includegraphics[width=.55\linewidth]{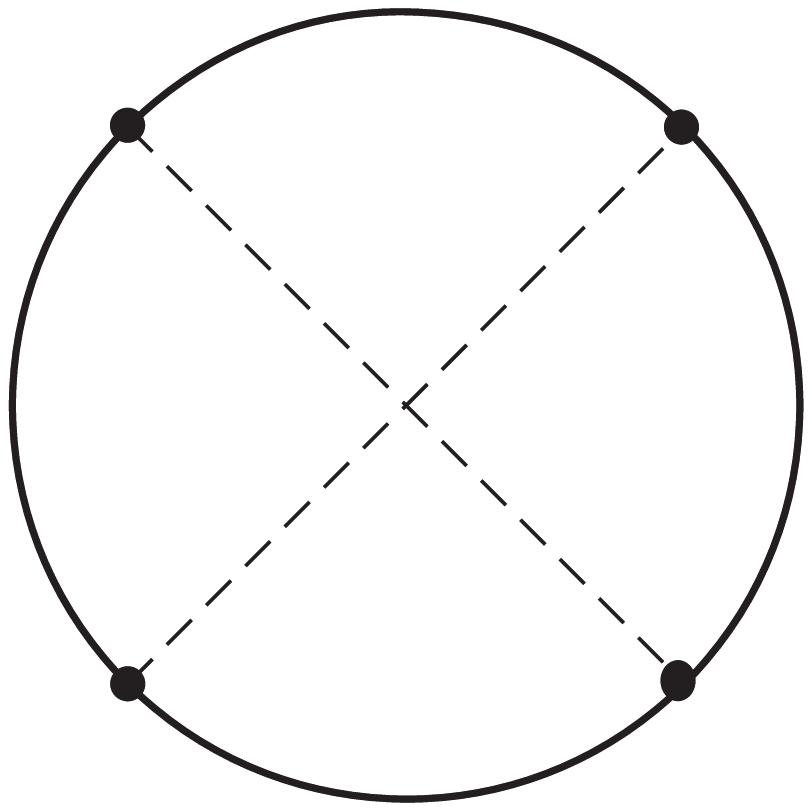}
\caption{The ``$X$-graph"}
\label{fig_X-graph}
\end{center}
\end{minipage}
\hskip 0.4cm
\begin{minipage}{.45\linewidth}
\begin{center}
\includegraphics[width=.55\linewidth]{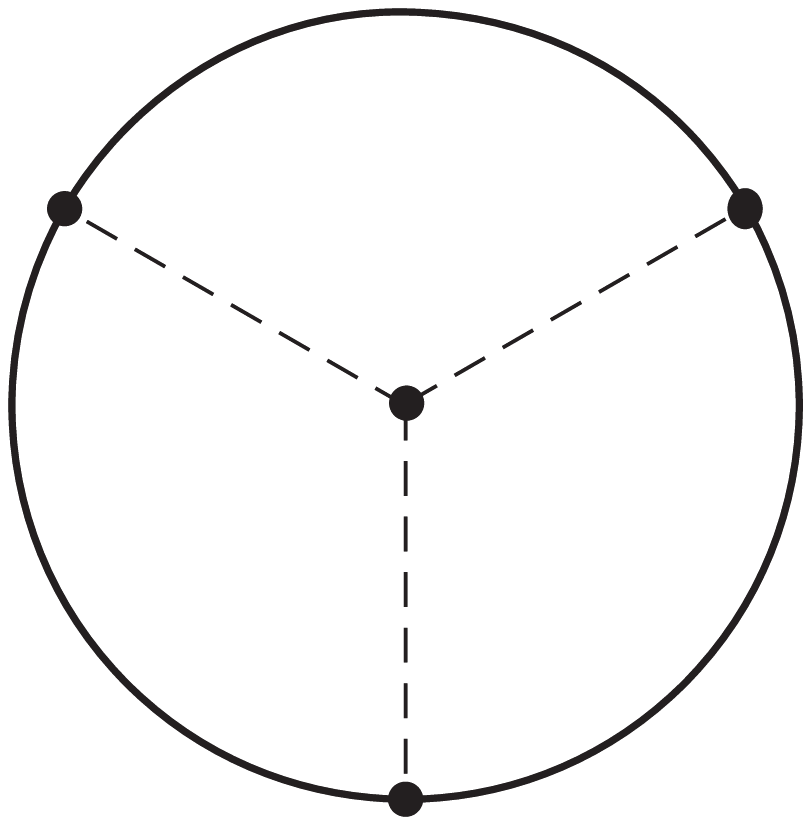}
\caption{The ``$Y$-graph"}
\label{fig_Y-graph}
\end{center}
\end{minipage}
\end{center}
\end{figure}

\begin{theorem} {\rm (\cite{Bo-Ta}, \cite{Li-Wa})} 
There holds 
$$
\frac14 GI_X(K)-\frac13 GI_Y(K)+\frac1{24}=v_2(K),
$$
where $v_2(K)$ is the second coefficient of the Conway polynomial. 
\end{theorem}

\subsection{Expression in terms of $1$-form}\label{ss_expression}

In order to generalize $GI_Y(K)$ in a conformal geometric way in Subsection \ref{ss_conf_inv_Y}, we give another expression of $GI_Y(K)$. 

For a vector $v\in T_y\mathbb{R}^3$ 
we define a vector field $X_{G}(y;v)$ on $\mathbb{R}^3\setminus\{y\}$ by 
$$
X_{G}(y;v)(x)=\frac1{4\pi}\cdot\frac{x-y}{|x-y|^3}\times v 
=\frac1{4\pi}\, v\times\nabla\left(\frac1{|x-y|}\right)
\hskip 0.5cm (x\in \mathbb{R}^3\setminus\{y\}), 
$$
and a $1$-form $\lambda_{G} (y;v)$ on $\mathbb{R}^3\setminus\{y\}$ by 
\[\lambda_{G} (y;v)(u)=u\cdot X_{G}(y;v)(x)
=\frac{\det(x-y,v,u)}{4\pi|x-y|^3}
\hskip 0.5cm (u\in T_x\mathbb{R}^3).\] 

Let $v_y$ denote a unit tangent vector to $K$ at $y$. 
Then $GI_Y(K)$ can be expressed in terms of this $1$-form, 
or equivalently, by the vector field as follows (\cite{Li-Wa}):
\begin{equation}\label{f_GI_Y}
\begin{array}{l}
GI_Y(K) =\displaystyle \int_{U_Y}
dy_1dy_2dy_3\,\lambda_{G}(y_1\,;v_{y_1})\w\lambda_{G}(y_2\,;v_{y_2})\w\lambda_{G}(y_3\,;v_{y_3})\, \\[2mm]
=-\displaystyle \int_{U_Y}
\det \left(X_{G}(y_1;v_{y_1})(x), X_{G}(y_2;v_{y_2})(x), X_{G}(y_3;v_{y_3})(x)\right)
\,dy_1dy_2dy_3d_{\rm vol}x. 
\end{array}
\end{equation}

\begin{definition}\label{def_lambda_K} \rm 
Let $K=f(S^1)$ be a knot. 
Define a $1$-form $\lambda_K$ on the complement of the knot $K$ by 
$$
\lambda_K(u)=\displaystyle{\int_K\lambda_{G} (y;v_y)(u)dy}
=\frac1{4\pi}\displaystyle{\int_K\frac{\det(x-y,v_y,u)}{|x-y|^3}\,dy}
$$
for $u\in T_x\mathbb{R}^3$ $(x\not\in K)$, where $v_y$ denotes the unit tangent vector to $K$ at $y$. 
\end{definition}
Although $\lambda_{G} (y;v_y)$ is not closed, 
\begin{lemma}\label{lambda_K=closed}
The $1$-form $\lambda_K$ is closed. 
\end{lemma}
\begin{proof} 
Suppose the knot $K$ is given by $K=f(S^1)$. 
Then 
\[4\pi\lambda_K(u)=u\cdot\int_{S^1}\frac{x-f(s)}{|x-f(s)|^3}\times f^{\p}(s)\,ds.\]
Therefore, the $2$-form $d\lambda_K$ vanishes if and only if 
\[\nabla\times\left(\int_{S^1}\frac{x-f(s)}{|x-f(s)|^3}\times f^{\p}(s)\,ds\right)=\vect 0, \]
which is the consequence of 
\[\begin{array}{rcl}
\displaystyle \nabla\times\left(\frac{x-f(s)}{|x-f(s)|^3}\times f^{\p}(s)\right)
&=&\displaystyle -(f^{\p}(s)\cdot\nabla)\left(\frac{x-f(s)}{|x-f(s)|^3}\right) \\[4mm]
&&\displaystyle  +\nabla^2\left(\frac1{|x-f(s)|}\right)f^{\p}(s)\\[4mm]
&=&\displaystyle \frac{f^{\p}(s)}{|x-f(s)|^3}-3\,\frac{f^{\p}(s)\cdot(x-f(s))}{|x-f(s)|^5}(x-f(s))\\[4mm]
&=&\displaystyle -\frac{d}{ds}\left(\frac{x-f(s)}{|x-f(s)|^3}\right).
\end{array}\]
\end{proof}

When $K^{\prime}$ is a knot in $\mathbb{R}^3\setminus K$ then 
$$\int_{K^{\prime}}\lambda_K=\frac1{4\pi}\displaystyle{\iint_{K\times K^{\prime}}\frac{\det(x-y,v_y,v_x)}{|x-y|^3}\,dxdy}$$
is the linking number of $K$ and $K^{\prime}$. 

\section{Preliminaries from conformal geometry} 

\begin{definition}\label{def_conf_angle} \rm (Doyle and Schramm) 
Let $C(x,y,y)$ $(x\ne y)$ be the circle which is tangent to $K$ at $y$ that passes through $x$, and let $v_x$ be the unit tangent vector to $K$ at $x$. 
Define the conformal angle $\theta_K(x,y)$ $(0\le\theta_K(x,y)\le\pi)$ by the angle between $C(x,y,y)$ and $v_x$ at $x$. 
\end{definition}

\begin{lemma}\label{order_conf_angle}{\rm (\cite{La-OH})} 
The conformal angle $\theta _K(x,y)$ is of the order of $|x-y|^2$ near the diagonal. 
To be precise, 
\[\theta _K(x,y)=\frac{\sqrt{{\kappa^{\p}}^2+\kappa^2\tau^2\,}}6\,|x-y|^2+O(|x-y|^3)\]
near the diagonal, where $s, \kappa, \tau$ denote the arc-length, curvature, and torsion of $K$ respectively, and $\kappa^{\p}$ means $\displaystyle \frac{d\kappa}{ds}$. 
\end{lemma}

\begin{definition} \rm 
Let $T$ be a M\"obius transformation of $\mathbb{R}^3\cup\{\infty\}$. 
Define $|T^{\p}(r)|$ for $p\in\mathbb{R}^3$ by 
$$|T^{\p}(r)| =\left|\det d\,T(p)\right| ^{\frac{\,1\,}{6}}.$$
\end{definition}
\begin{lemma}\label{lem_formulae_Moeb} 
\begin{enumerate}
\item If $I_0(p)$ is an inversion in a sphere of radius $r$ with center the origin then 
\begin{equation}\label{f_I_r^p}
|I_0(p)^{\p}(p)|=\frac r{\,|p|\,}. 
\end{equation}
\item We have 
\begin{equation}\label{f_|T(p)-T(q)|}
|T(p)-T(q)|=\left|T^{\p}(p)\right|\left|T^{\p}(q)\right||p-q|
\end{equation}
for a pair of points $p,q$ in $\mathbb{R}^3$ and 
\begin{equation}\label{f_T(gamma)'}
|T_{\ast}v|=|T^{\p}(p)|^2|v|
\end{equation}
for a vector $v$ in $T_p\mathbb R^3$. 
\end{enumerate}
\end{lemma}
\begin{corollary}\label{conf_inv_2-form}
The $2$-form $\displaystyle \frac{dxdy}{|x-y|^2}$ on $K\times K\setminus\Delta$ is conformally invariant. 
In other words, if $T$ is a M\"obius transformation and $\tilde x$ and $\tilde y$ denote $T(x)$ and $T(y)$ then 
\begin{equation}\label{f_conf_inv_|inf_cr|}
T^{\ast}\left(\frac{d\tilde xd\tilde y}{|\tilde x-\tilde y|^2}\right)=\frac{dxdy}{|x-y|^2}.
\end{equation}
\end{corollary}

\section{Conformally invariant $1$-form via vector field}

We shall construct a vector field on $\mathbb R^3\setminus\{y\}$ for a vector in $T_y\mathbb{R}^3$ in a conformally invariant manner. 
It gives a confomally invariant $1$-form on $\mathbb R^3\setminus\{y\}$. 
We can define confomally invariant functionals associated with the $Y$-graph in terms of this $1$-form. 

We make use of tangent circles at a given vector in $T\mathbb{R}^3$ since a M\"obius transformation maps a circle into a circle and tangent curves into tangent curves. 

\begin{definition}\label{def_conf_transport} \rm 
Let $x$ and $y$ be a pair of points $(x\ne y)$ and $v$ a vector in $T_y\mathbb{R}^3$. 
Let $C(y;v,x)$ denote the circle through $x$ and $y$ which is tangent to $v$ at $y$ whose orientation is given by $v$. 

\begin{figure}[htbp]
\begin{center}
\includegraphics[width=.3\linewidth]{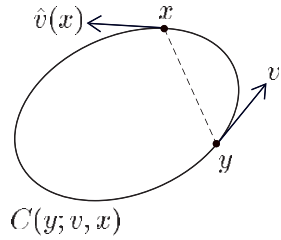}
\caption{$v$ and $\hat v$}
\label{conf transp}
\end{center}
\end{figure}

\smallskip
(1)
Let $\hat v(x)$ be a tangent vector to $C(y;v,x)$ at $x$ with the same norm as $v$. 
As $\hat v(x)$ is symmetric to $v$ in the line joining $x$ and $y$, it is given by 
\begin{equation}\hat v(x)=\displaystyle 
\left\{2\left(v\cdot \frac{x-y}{|x-y|}\right)\frac{x-y}{|x-y|}-v\right\}.
\end{equation}
We call the correspondence $v\mapsto \hat v(x)$ the {\em unit conformal transportation} from $y$ to $x$ (\cite{La-OH}). 

\smallskip
(2) Put
$$\tilde v(x)=\frac{\hat v(x)}{\,|x-y|^2\,}
=\displaystyle \frac {1}{\,|x-y|^2\,}
\left\{2\left(v\cdot \frac{x-y}{|x-y|}\right)\frac{x-y}{|x-y|}-v\right\}
$$
for $x,y\in\mathbb{R}^3$ $(x\ne y)$ and $v\in T_y\mathbb{R}^3$, 
and call it the {\em vector $v$ conformally transported to $x$}. 
We call the correspondence 
$$T_y\mathbb{R}^3\ni v\mapsto \tilde v(x)\in T_x\mathbb{R}^3$$
the {\em conformal transportation} from $y$ to $x$.
\end{definition}

Once a vector field is given, a $1$-form can be defined in the same way as before. 

\begin{definition}\label{def_widetilde_omega}
Let $v$ be a venctor in $T_y\mathbb{R}^3$. 
Define a $1$-form $\widetilde\omega=\widetilde\omega(y\,;v)$ on $\mathbb{R}^3\setminus\{y\}$ by 
$$\widetilde\omega(u)=u\cdot \tilde v(x) \hskip 0.5cm (x\in \mathbb{R}^3\setminus \{y\}, \, u\in T_x\mathbb{R}^3).$$
\end{definition}

\begin{remark} \rm 
The $1$-form $\hat\omega=\hat\omega(y\,;v)$ on $\mathbb{R}^3\setminus\{y\}$, where $v$ is a unit vector in $T_y\mathbb R^3$, defined by 
$$\hat\omega(u)=u\cdot \hat v(x) \hskip 0.5cm (x\in \mathbb{R}^3\setminus \{y\}, \, u\in T_x\mathbb{R}^3)$$
was used by H\'elein to show the isoperimetric inequality (\cite{He}). 
\end{remark}
Then the vector field and the $1$-form can be expressed 
in terms of an inversion in a $2$-sphere. 
\begin{lemma}\label{lemma}
Let $I_p$ $(p\in\mathbb{R}^3)$ denote an inversion in the $2$-sphere with center $p$ and raidus $1$. 
Let $v$ be a vector in $T_y\mathbb R^3$ and $x$ a point in $\mathbb R\setminus\{y\}$. 

\smallskip
{\rm (1)} $\tilde v(x)=-{I_x}_{\ast}(v)$ as a vector. 

\smallskip
{\rm (2)} $\widetilde\omega(y\,;v)(u)=-v\cdot ({I_y}_{\ast}u)$ for $u\in T_x\mathbb{R}^3$. 
\end{lemma}
\begin{proof}
%
(1) Since $v$ is tangent to the circle $C(y;v,x)$, the vector ${I_x}_{\ast}v$ is tangent to ${I_x}(C(y;v,x))$. 
Since the tangent vector to $C(y;v,x)$ at $x$ is equal to a positive multiple of $\tilde v(x)$, ${I_x}(C(y;v,x))$ is a line whose direction vector is equal to a positive multiple of $-\tilde v(x)$. 
Therefore, ${I_x}_{\ast}v$ is a positive multiple of $-\tilde v(x)$. 
Then the conclusion comes from 
$\displaystyle |{I_x}_{\ast}v|=\nekofrac{|v|}{|x-y|^2}=|\tilde v(x)|,$ 
which is implied by Lemma \ref{lem_formulae_Moeb}. 

\smallskip
(2) The image ${I_y}(C(y;v,x))$ is the line through $I_y(x)$ whose direction vector is equal to a positive multiple of $-v$. 
Since $\tilde v(x)$ is tangent to $C(y;v,x)$, it follows that ${I_y}_{\ast}\tilde v(x)$ is a positive multiple of $-v$. 
Since $\angle {I_y}_{\ast}$ is a conformal mapping, we have the following equality between angles: 
\[\angle {I_y}_{\ast}u\cdot (-v)=\angle {I_y}_{\ast}u\cdot {I_y}_{\ast}\tilde v(x)=\angle u\cdot \tilde v(x).\] 
Since $|{I_y}_{\ast}u|=\nekofrac{|u|}{|x-y|^2}$ by Lemma \ref{lem_formulae_Moeb}, we have 
\[\widetilde\omega(y\,;v)(u)=u\cdot \tilde v(x)=({I_y}_{\ast}u)\cdot (-v).\]
\end{proof}

As a corollary we have 

\begin{proposition} \label{omega=exact} 
Define a map $\psi:\mathbb{R}^3\setminus \{y\}\to\mathbb{R}$ by 
\[\psi (x)=-I_y(x)\cdot v.\] 
Then $\widetilde\omega(y\,;v)=d\psi$, i.e. the $1$-form $\widetilde\omega(y\,;v)$ is exact. 
\end{proposition}
\begin{proof} Suppose $x\in \mathbb R\setminus\{y\}$ and $u\in T_x\mathbb R^3$. 
Then Lemma \ref{lemma} (2) implies 
\[d\psi (u)=-({I_y}_{\ast}u)\cdot v=\widetilde\omega(y\,;v)(u).\] 
\end{proof}

We remark that the $1$-form $\lambda_{G} (y;v_y)$ is not necessarily closed. 

\begin{proposition} 
The $\widetilde\omega(y\,;v)$ is conformally invariant, i.e. 
\[\widetilde\omega(T(y)\,;T_{\ast}v)(T_{\ast}u)=\widetilde\omega(y\,;v)(u)\] 
for any $u\in T_x\mathbb{R}^3$ $(x\in\mathbb R^3\setminus\{y\})$ and M\"obius transformation $T$ of $\mathbb{R}^3\cup\{\infty\}$.
\end{proposition}

\begin{proof} 
Lemma \ref{lem_formulae_Moeb} implies 
\begin{eqnarray}
|T_{\ast}u|&=&|T^{\p}(x)|^2|u|, \label{f_T_ast_u}\\
|T_{\ast}v|&=&|T^{\p}(y)|^2|v|.\nonumber\\
\displaystyle \left|\widetilde{T_{\ast}v}(T(x))\right|&=&\displaystyle \frac{|T_{\ast}v|}{|T(y)-T(x)|^2}=\frac{|v|}{|T^{\p}(x)|^2|x-y|^2}, \label{f_tilde_T_T_x}
\end{eqnarray}
where $\widetilde{T_{\ast}v}(T(x))$ denotes the vector $T_{\ast}v$ conformally transported to $T(x)$ (Definition \ref{def_conf_transport}). 

We have 
\begin{equation}\label{f_angle}
\angle T_{\ast}u\cdot \widetilde{T_{\ast}v}(T(x))=\angle u\cdot \tilde v(x).\end{equation}
This is because the right hand side is equal to the angle between $u$ and $C(y;v,x)$ at $x$ as $\tilde v(x)$ is tangent to the circle $C(y;v,x)$ at $x$, whereas the left hand side is equal to the angle between $T_{\ast}u$ and $T(C(y;v,x))$ at $T(x)$ as $\widetilde{T_{\ast}v}(T(x))$ is tangent to the circle $T(C(y;v,x))$ at $T(x)$. 

Putting (\ref{f_T_ast_u}), (\ref{f_tilde_T_T_x}), and (\ref{f_angle}) together we have 
\[\widetilde\omega(T(y)\,;T_{\ast}v)(T_{\ast}u)=(T_{\ast}u)\cdot \widetilde{T_{\ast}v}(T(x))=u\cdot \tilde v(x)=\widetilde\omega(y\,;v)(u).\]
\end{proof}

\medskip
Let us play the same game as in the Subsection \ref{Lin-Wang} by substituting the conformally invariant $1$-form $\widetilde\omega$ for the $1$-form $\lambda_G$ coming from the Gauss formula for the linking number. 
But unfortunately, it turns out that this attempt is not successful. 

A similar construction of $\lambda_K$ from $\lambda_{G} (y;v_y)$ gives a trivial $2$-form for $\widetilde\omega(y\,;v_y)$. 

\begin{proposition}
Let $\widetilde\omega_K$ be a $1$-form on $\mathbb{R}^3\setminus K$ defined by 
$$
\widetilde\omega_K(u)=\displaystyle{\int_K\widetilde\omega(y;v_y)(u)dy}
=\displaystyle \int_K u\cdot \tilde v_y(x)\,dy
=\displaystyle u\cdot \int_{S^1}\widetilde{f^{\prime}(s)}(x)\,ds
$$
for $u\in T_x\mathbb{R}^3$, where $v_y$ is the unit tangent vector to $K$ at $y$. 
Then $\widetilde\omega_K$ vanishes for any knot $K$. 
\end{proposition}

\begin{proof}
Lemma \ref{lemma} indicates that 
\begin{equation}
\int_{S^1}\widetilde{f^{\prime}(s)}(x)\,ds
=-\int_{S^1}{I_x}_{\ast}\left(f^{\prime}(s)\right)ds
=-\int_{S^1}\frac{d}{ds}\left(I_x\circ f(s)\right)ds
=0. 
\label{eq_widetilde_f}
\end{equation}
\end{proof}

\section{Conformally invariant $Y$-energy}\label{ss_conf_inv_Y}

Using the conformaly invariant $1$-form given in the previous section 
we can consider the same construction of $GI_Y$ from $\lambda_K$ (\ref{f_GI_Y}). 
But unfortunately, it gives a trivial functional. 

Let $v_{y_i}$ ($i=1,2,3)$ be a unit tangent vector at $y_i$ and $\hat v_{y_i}=\hat v_{y_i}(x)$ be the image of $v_{y_i}$ under the unit conformal transportation from $y_i$ to $x$ (Definition \ref{def_conf_transport}). 
Let $\widetilde\omega(y\,;v)$ be the $1$-form given in Definition \ref{def_widetilde_omega}. 
Let $U_Y$ be the domain in $K^3\times\mathbb R^3$ given by (\ref{def_U_Y}). 

\begin{proposition} 
Put 
\begin{eqnarray}\label{f_Y-energy}
E_Y^{\circ}(K)&=&\displaystyle \int_{U_Y}
dy_1dy_2dy_3\,\widetilde\omega(y_1\,;v_{y_1})\w\widetilde\omega(y_2\,;v_{y_2})\w\widetilde\omega(y_3\,;v_{y_3})\nonumber \\[2mm]
&=&\displaystyle \int_{U_Y}\frac{\det(\hat v_{y_1}(x), \hat v_{y_2}(x), \hat v_{y_3}(x))}{|y_1-x|^2|y_2-x|^2|y_3-x|^2}
\,dy_1dy_2dy_3d_{\rm vol}x. \hskip 0.5cm {}
\end{eqnarray}
Then $E_Y^{\circ}(K)=0$ for any knot $K$. 
\end{proposition}

It is enough to show the following. 

\begin{lemma} 
Suppose $y_1\ne y_2\ne y_3\ne y_1$. 
Then 
\[\int_{\mathbb{R}^3\setminus \{y_1, y_2, y_3\}}
\widetilde\omega(y_1\,;v_{y_1})\w\widetilde\omega(y_2\,;v_{y_2})\w\widetilde\omega(y_3\,;v_{y_3})
=0.\]
\end{lemma}

\begin{proof} 
Define a map $\psi_i:\mathbb{R}^3\setminus \{y_i\}\to\mathbb{R}$ $(i=1,2,3)$ by 
\[\psi_i (x)=-I_{y_i}(x)\cdot v_{y_i},\]
where $I_{y_i}$ is an inversion in a sphere with radius $1$ and center $y_i$. 
Proposition \ref{omega=exact} implies that 
\[\omega(y_i\,;v_{y_i})=d\psi_i.\]
Let $\Si_i$ $(i=1,2,3)$ be a sphere with radius $\e$ $(\e\ll 1)$ and center $y_i$, and $\Si_0$ a sphere with radius $R$ $(R\gg 1)$ and center the origin. 
Let $\Omega_{\e,R}$ be the domain in $\mathbb R^3$ bounded by $\Si_0$, $\Si_1$, $\Si_2$, and $\Si_3$. 
Then 
\[\begin{array}{l}
\displaystyle \int_{\mathbb{R}^3\setminus \{y_1, y_2, y_3\}}
\widetilde\omega(y_1\,;v_{y_1})\w\widetilde\omega(y_2\,;v_{y_2})\w\widetilde\omega(y_3\,;v_{y_3})\\[4mm]
=\displaystyle \lim_{\e\to+0, R\to+\infty}\int_{\Omega_{\e,R}}d(\psi_1d\psi_2\w d\psi_3)\\[4mm]
=\displaystyle \lim_{\e\to+0, R\to+\infty}\left(\int_{\Si_0}-\int_{\Si_1}-\int_{\Si_2}-\int_{\Si_3}\right)\psi_1d\psi_2\w d\psi_3.
\end{array}\]

\smallskip
(i) On $\Si_0$ we have $|\psi_1|=O(\frac1R)$ and $|d\psi_2\w d\psi_3|=O(\frac1{R^4})$ whereas $\textrm{Area}(\Si_0)=O(R^2)$. 
Therefore 
\[\lim_{R\to+\infty}\int_{\Si_0}\psi_1d\psi_2\w d\psi_3=0.\]

\smallskip
(ii) On $\Si_2$ we have $\psi_1=C+O(\e)$ for some constant $C$. 
Therefore 
\[\begin{array}{rcl}
\displaystyle \lim_{\e\to+0}\int_{\Si_2}\psi_1d\psi_2\w d\psi_3
&=&\displaystyle C\lim_{\e\to+0}\int_{\Si_2}d\psi_2\w d\psi_3\\[4mm]
&=&\displaystyle C\lim_{\e\to+0}\int_{\Si_2}d(\psi_2\w d\psi_3)\\[4mm]
&=&0.
\end{array}\]
The same argument works for $\Si_3$. 

\smallskip
(iii) On $\Si_1$ we have $|\psi_1|=O(\frac1{\e})$ and 
\[d\psi_2\w d\psi_3=\omega^{\p}+O(\e)\]
for some constant $2$-form $\omega^{\p}$, whereas $\textrm{Area}(\Si_1)=O(\e^2)$. 
Therefore 
\[\lim_{\e\to+0}\int_{\Si_1}\psi_1d\psi_2\w d\psi_3=0.\]
\end{proof}

\begin{proposition}
Put 
\begin{equation}\label{f_abs-Y-energy}
AE_Y^{\circ}(K)
=\displaystyle \int_{U_Y}\frac{\left|\det(\hat v_{y_1}(x), \hat v_{y_2}(x), \hat v_{y_3}(x))\right|}{|y_1-x|^2|y_2-x|^2|y_3-x|^2}
\,dy_1dy_2dy_3d_{\rm vol}x. 
\end{equation}
Then $AE_Y^{\circ}$ diverges for any knot $K$. 
\end{proposition}

\begin{proof}
The integrand blows up near the diagonal of $K^3\times \mathbb{R}^3$, 
which is stratified into strata of different (co)dimensions. 
The contribution of the neibhborhood of 
$${\mathcal N}_2=\{(y_1,y_2,y_3,x)\in K^3\times \mathbb{R}^3\,|\, x=y_1=y_2\ne y_3\}$$ 
makes the integral diverge. 

Let us fix $y_1$ and $y_3$ $(y_1\ne y_3)$. 
We may assume without loss of generality that the knot $K$ is parametrized by the arc-length as $K=\gamma([0,1])$ and that $y_1=\gamma(0)=\vect 0$. 
Let $\Omega(\e)$ be a domain given by 
\[\Omega(\e)=\left\{(y_2, x)\in K\times \mathbb R^3:
y_2=\gamma(s), x\ne y_2, \frac{\,\e\,}2<s, |x|\le \e
\right\}.\]
Since $v_{y_3}(x)=v^{\p}+O(\e)$ on $\Omega(\e)$ for some constant vector $v^{\p}$, we have the similarity between the integrands on $\Omega(\e)$ and $\Omega(\frac{\e}2)$: 
\[\frac{\big|\det(\hat v_{0}(\frac{x}2), \hat v_{\gamma(\frac{s}2)}(\frac{x}2), \hat v_{y_3}(\frac{x}2))\big|}{|\frac{x}2|^2|\gamma(\frac{s}2)-\frac{x}2|^2|y_3-\frac{x}2|^2}
=2^4\frac{\left|\det(\hat v_{0}(x), \hat v_{\gamma(s)}(x), \hat v_{y_3}(x))\right|}{|x|^2|\gamma(s)-x|^2|y_3-x|^2}+O(\e^{-3}).\]
Since the volume of $\Omega(\e)$ is of order $\e^4$ we have 
\[\begin{array}{l}
\displaystyle 
\int_{\Omega(\frac{\e}2)}\frac{\left|\det(\hat v_{0}(x), \hat v_{y_2}(x), \hat v_{y_3}(x))\right|}{|x|^2|y_2-x|^2|y_3-x|^2}\,dy_2d_{\rm vol}x\\[4mm]
=\displaystyle \int_{\Omega(\e)}\frac{\left|\det(\hat v_{0}(x), \hat v_{y_2}(x), \hat v_{y_3}(x))\right|}{|x|^2|y_2-x|^2|y_3-x|^2}\,dy_2d_{\rm vol}x+O(\e).
\end{array}\]
It follows that there is a positive constant $C^{\p}$ such that 
\[\int_{\Omega(2^{-n})}\frac{\left|\det(\hat v_{0}(x), \hat v_{y_2}(x), \hat v_{y_3}(x))\right|}{|x|^2|y_2-x|^2|y_3-x|^2}\,dy_2d_{\rm vol}x\ge C^{\p}.\]
for all $n\in\mathbb N$. 
Therefore
\[\begin{array}{rcl}
AE_Y^{\circ}(K)&\ge&\displaystyle \int\left(\sum_{n=1}^{\infty}\int_{\Omega(2^{-n})}\frac{\left|\det(\hat v_{y_1}(x), \hat v_{y_2}(x), \hat v_{y_3}(x))\right|}{|y_1-x|^2|y_2-x|^2|y_3-x|^2}\,dy_2d_{\rm vol}x\right)dy_1dy_3\\[4mm]
&=&\infty.
\end{array}\]
\end{proof}

\begin{remark} \rm 
There is a well-defined non-trivial conformally invariant ``trilocal'' functional on the space of knots (\cite{La-OH2}). 
Although it is not associated to a Y-diagram, the integrand involves a triplet $(x, v_x), (y, v_y)$, and $(z, v_z)$. 
\end{remark}

\bigskip
\noindent
{\large{\bf Acknowledgement}}
The author thanks R\'emi Langevin deeply for many helpful suggestions.

Department of Mathematics, Tokyo Metropolitan University, 

1-1 Minami-Ohsawa, Hachiouji-Shi, Tokyo 192-0397, JAPAN. 

ohara@comp.metro-u.ac.jp

\end{document}